\documentclass[12pt]{article}
\usepackage{latexsym, amsmath, amsfonts}

\newcommand{\R}{{\mathbb R}}
\newcommand{\reel}{{\mathbb R}}                         
\newcommand{\entier}{{\mathbb N}}                      

\newcommand{\dis}{\displaystyle}
\newcommand{\E}{{\mathbb E}}
\newcommand{\Prob}{{\mathbb P}}

\newcommand{\B}{{\cal B}}
\newcommand{\lan}{\langle}
\newcommand{\ran}{\rangle}
\newcommand{\gr}{\kappa}
\newcommand{\Id}{\rm Id}

\newcommand{\al}{\alpha}
\newcommand{\eps}{\varepsilon}
\newcommand{\la}{\lambda}
\newcommand{\de}{\delta}

\newcommand{\pr}[2]{\langle {#1} , {#2} \rangle}
\newcommand{\norm}[1]{\left \| #1 \right \|}

\newtheorem{theorem}{Theorem}
\newtheorem{corollary}{Corollary}
\newtheorem{lemma}{Lemma}
\newtheorem{proposition}{Proposition}
\newenvironment{thm}[2][]
               {\pagebreak[3] \medskip {\bf \noindent Theorem}
                               {#1} {\it #2}
                             }{}
\title{$L_p$ moments of random vectors via majorizing measures}
\author{Olivier Gu\'edon,
Mark Rudelson \thanks{Research was supported in part by  NSF grant
DMS-024380.}}
\date{}

\begin{document}
\maketitle
\begin{abstract}
For a  random vector $X$ in $\R^n$, we obtain bounds on the size
of a sample, for which the empirical $p$-th moments of linear
functionals are close to the exact  ones  uniformly on a convex
body $K \subset \R^n$. We prove an estimate for a  general random
vector and apply it to several problems arising  in geometric
functional analysis. In particular, we find a short Lewis type
decomposition for any finite dimensional subspace of $L_p$. We
also prove that for an isotropic $\log$-concave random
vector, we  only need $\lfloor n^{p/2} \log n \rfloor$ sample
points so that the empirical $p$-th moments of the linear
functionals are almost isometrically the same as the exact ones.
We obtain a concentration estimate for the empirical moments. The
main ingredient of the proof is the construction of an appropriate
majorizing measure to bound a certain Gaussian process.
\end{abstract}

\section{Introduction}
In many problems of geometric functional analysis it is necessary
to approximate a given random vector by an empirical sample. More
precisely, given a random vector $X \in \R^n$, we want to find the
smallest number $m$ such that the properties of $X$ can be
recovered from the empirical measure $1/m \sum_{j=1}^m
\delta_{X_j}$, constructed with independent copies $X_1, \ldots,
X_m$ of the vector $X$. In particular, for $p \ge 2$ and for $y
\in \R^n$, we want to approximate the moments $\E |\lan X,y \ran
|^p$ by the empirical averages $1/m \sum_{j=1}^m |\lan X_j,y \ran
|^p$ with high probability. Moreover, we require this
approximation to be uniform over $y$ belonging to some convex
symmetric set in $\R^n$. A problem of this type was considered in
\cite{BouLinMil}. Formulated in analytic language, it asks about
finding the smallest $m$ and a set of points $x_1, \ldots x_m \in
X$ such that for any function $f$ from an $n$-dimensional function
space $F \subset L_1(X,\mu)$,
$$
  (1-\eps)\|f\|_1
  \le \frac{1}{m} \sum_{j=1}^m |f(x_j)|
  \le (1-\eps)\|f\|_1.
$$
Another example of such problems originates in Computer Science.
The probabilistic algorithm for estimating the volume of an
$n$-dimensional convex body, constructed by Kannan, Lov\'asz, and
Simonovits \cite{KLS} required to bring the body to a nearly
isotropic position as a preliminary step. To this end, one has to
sample $m$ random points $x_1, \ldots, x_m$ in the body $L$ so that
the empirical isotropy tensor will be close to the exact one,
namely
\begin{equation}  \label{inertia}
  \norm{\frac{1}{m} \sum_{j=1}^m x_j \otimes x_j
  - \frac{1}{\text{vol}(L)} \int_L x \otimes x \, dx}
  < \eps.
\end{equation}
This problem was attacked with different probabilistic techniques.
The original estimate of \cite{KLS} was significantly improved by
Bourgain \cite{Bou}. Using the decoupling method he proved that
$m=C(\eps)n \log^3 n$ vectors $x_1, \ldots x_m$ uniformly
distributed in the body $L$ satisfy \eqref{inertia} with high
probability. This estimate was farther improved to $(Cn/\eps^2)
\cdot \log^2 (Cn/\eps^2)$ in \cite{Rud2}, \cite{Rud3}. The proof
in \cite{Rud2} used majorizing measures, while the later proof in
\cite{Rud3} was based on the non-commutative Khinchine
inequality.

These problems were put into a general framework by Giannopoulos
and Milman \cite{GiaMil}, who related them to the concentration
properties of a random vector. Let $\al>0$ and let $\nu$ be a
probability measure on $(X,\Omega)$. For a function $f: X \to \R$
define the $\psi_{\al}$-norm by
\[
  \norm{f}_{\psi_\al}
  = \inf \{ \la>0 \mid \int_X \exp (|f|/\la)^{\al} d\nu \le 2 \}.
\]
Chebychev's inequality shows that the functions with bounded
$\psi_\al$-norm are strongly concentrated, namely $\nu\{x \mid
|f(x)|> \la t\} \le C \exp(-t^\al)$. Let $\mu$ be a Borel measure
in $\R^n$. It is called isotropic if
\[
  \int_{\R^n} x \otimes x \, d \mu(x) = \Id,
\]
where $\Id$ is the identity operator in $\R^n$. Note that this
normalization is consistent with the one used in \cite{KLS, Rud2,
Rud3}. The normalization  used in \cite{MilPaj, GiaMil} differs
from it by the multiplicative coefficient $L_{\mu}^2$, where
$L_{\mu}$ is the {\em isotropic constant} of $\mu$ (see
\cite{MilPaj}).

The paper \cite{GiaMil} considers isotropic measures which satisfy
the $\psi_\al$-condition for scalar products:
\[
  \norm{\pr{\cdot}{y}}_{\psi_\al} \le C
\]
for all $y \in S^{n-1}$. Here and below $C,c, \ldots$ denote
absolute constants, whose value may change at each occurrence.

Note that by Borell's lemma, any log-concave measure in $\R^n$
satisfies the $\psi_1$-condition \cite{MilSch}, \cite{MilPaj}. Let
$p \ge 1$ and let $\mu$ be an isotropic log-concave measure
satisfying the $\psi_\al$ condition for scalar products with some
$\al \in [1,2]$. The central result of \cite{GiaMil} provides an
estimate for the minimal size of a set of independent random
vectors $X_1, \ldots, X_m$ distributed according to the measure
$\mu$ such that the empirical $p$-moments satisfy the inequality
\begin{equation}   \label{GM result}
  \Gamma_1(p)
  \le \left ( \frac{1}{m} \sum_{j=1}^m |\pr{X}{y}|^p
      \right )^{1/p}
  \le \Gamma_2(p),
  \quad
  \forall y \in S^{n-1}.
\end{equation}
The $\psi_\al$-condition implies that the
$L_p(\mu)$ and $L_2(\mu)$-norms of the function $f_y(x)=\pr{x}{y}$
are equivalent. Thus the  inequality \eqref{GM result} means that
the empirical $p$-moment of $f_y$ is equivalent to the real
$p$-moment up to a constant coefficient.

In the present paper we use a different approach to this problem
based on the majorizing measure technique developed by Talagrand
\cite{Tal:genericchaining}. This approach lead to breakthrough
results in various problems in probabilistic combinatorics and
analysis (see \cite{Tal:genericchaining} and references therein).
In a similar context the majorizing measures were applied in
\cite{Rud4} to select small almost orthogonal submatrices of an
orthogonal matrix, and in \cite{Rud2} to prove the estimate
\eqref{inertia} with small $m$.

To state the results we have to introduce some notation. Let
$(\reel^n, \langle \cdot, \cdot \rangle)$ be a Euclidean space,
and let $|\cdot|_2$ be the associated Euclidean norm. For a
symmetric convex body $K$ in $\reel^n$, we denote by $\| \cdot
\|_K$ the norm, whose unit ball is $K$, and by $K^o = \{y\in \R^n
\mid \forall x \in K, \langle x,y \rangle \le 1 \}$ the polar of
$K$. We assume that the body $K$ has the modulus of convexity of
power type $q \ge 2$ (see Section \ref{generalrandomvector} for
the definition).
 Classical examples of
convex bodies satisfying this property are unit balls of finite
dimensional subspaces of $L_q$  \cite{Clarkson} or of
non-commutative $L_q$-spaces (like Schatten trace class matrices
\cite{Tomczak}).
We denote by $D$ the radius of the symmetric convex set $K$ i.e.
the smallest $D$ such that $K \subset D B_{2}^n$. For every $1 \le
q \le +\infty$, we define $q^*$ to be the conjugate of $q$, i.e.
$1/q + 1/q^* = 1$.

Given a random  vector $X$ in $\reel^n$, let $X_1, \ldots, X_m$ be
$m$ independent copies of $X$. Let $K \subset \R^n$ be a convex
symmetric body. Denote by
$$
V_{p}(K) = \sup_{y \in K} \bigg| \frac{1}{m} \sum_{j=1}^m |\lan
X_j, y \ran|^p \ - \ \E |\lan X, y \ran|^p \bigg|
$$
the maximal deviation of  the empirical $p$-moment of $X$ from the
exact one. We would like to bound $V_p(K)$ under minimal
assumptions on the body $K$ and random vector $X$. This will allow
us to choose the size of the sample $m$ for which this deviation
is small with high probability. Although the resulting statement
is pretty technical, it is applicable to a wide range of problems
arising in geometric functional analysis. We discuss some examples
in Sections \ref{subspace}, \ref{seclogconcave}.

To bound such  random process,  we must have some control of the
random variable $\max_{1 \le j \le m} |X_j|_2$. To this end we
introduce the parameter $\gr_{p,m}(X)$, which plays a key role
below
$$
\gr_{p,m}(X) = \bigg( \E \max_{1 \le j \le m} |X_j|_2^p
\bigg)^{1/p}.
$$
We prove  the following estimate for $V_p(K)$.
\begin{theorem}\label{thrandomintro}
Let $K \subset (\R^n, \langle \cdot, \cdot \rangle)$ be a
symmetric convex body of radius $D$. Assume that $K$ has modulus
of convexity of power type $q$ for some $q \ge 2$. Let $p \ge q$
and let $q^*$ be the conjugate of $q$.
\\
Let $X$ be a random vector in $\R^n$, and let $X_1, \ldots, X_m$
be independent copies of $X$. Assume that
$$
  C_{p,\lambda}  \frac{(\log m)^{2/q^*}}{ m} \dis (D \cdot \gr_{p,m}(X))^{p}
\le
 \de^2 \cdot \sup_{y \in K} \E|\lan X, y \ran|^p
$$
for some $\de<1$. Then
$$
\E V_p(K) \le 2 \de \cdot \sup_{y \in K} \E|\lan X, y \ran|^p.
$$
\end{theorem}
The constant $C_{p,\lambda}$ in Theorem 1  depends on $p$ and on
the parameter $\lambda$ in the definition of the modulus of
convexity of power type $q$ (see Section 2.1 for the definition).

Note that minimal assumptions on the vector $X$ are enough to
guarantee that $\E V_p(K)$ becomes small for large $m$. Indeed,
assume that the variable $|X|_2$ possesses a finite moment of
order $p + \eps$ for some positive $\eps$. Then
$$
\gr_{p,m}(X) \le \left(\E \sum_{j=1}^m |X_j|_2^{p+\eps}\right)^{1/p+\eps}
\le m^{1/p+\eps} (\E |X|_2^{p+\eps})^{1/p+\eps},
$$
so the quantity
$$
 \frac{(\log m)^{2/q^*}}{ m}  \cdot \gr_{p,m}^{p}(X)
$$
tends to $0$ when $m$ goes to $\infty$. Moreover, in most cases,
$\gr_{p,m}(X)$ may be bounded  by a simpler quantity:
\begin{equation}
\label{simplerMs}
\gr_{p,m}(X) \le \left( \E \sum_{j=1}^m |X_j|_2^{p}Ê\right)^{1/p}
\le e \left( \E |X|^s \right )^{1/s} =: e M_s,
\end{equation}
where $s=\max(p, \log m)$.

Theorem \ref{thrandomintro} improves the results of \cite{GiaMil}
in two ways. First, it contains an almost isometric approximation
of the $L_p$-moments of the random vector by empirical samples
(see Theorem \ref{logconcaveintro} below).
Second, the assumption on the norm of a random vector $X$ used in
Theorem \ref{thrandomintro} is weaker than the
$\psi_{\alpha}$-assumption on the scalar products, appearing in
\cite{GiaMil}. This allows to handle the situations, where the
$\psi_{\alpha}$-estimate does not hold (see e.g. approximate Lewis
decompositions, discussed in Section \ref{subspace}).

While Theorem \ref{thrandomintro} combined with Chebychev's
inequality provides a bound for $V_p(K)$, which holds with high
probability, it is often useful to have this probability
exponentially close to $1$. Using a measure concentration result
of Talagrand (\cite{LT} Theorem 6.21), we obtain such probability
estimate in Theorem \ref{concentration}.

We apply Theorem \ref{thrandomintro} to isotropic $\log$-concave
random vectors. This class  includes many naturally arising types
of random vectors, in particular a vector uniformly distributed in
an isotropic convex body (see Section \ref{seclogconcave} for
exact definitions). The empirical moments of $\log$-concave
vectors have been extensively studied in the last years
\cite{KLS}, \cite{Bou}, \cite{Rud3}, \cite{GiaMil}, \cite{GHT}.
 We will prove the following
\begin{theorem}\label{logconcaveintro}
For any $\eps \in (0,1)$ and $p \ge 2$ there exists $n_0(\eps, p)$
such that for any $n \ge n_0(\eps, p)$, the following holds:
let $X$ be a log-concave isotropic random vector in $\R^n$, let
 $X_1, \ldots, X_m$ be  independent copies of
$X$, if
$$
m = \lfloor C_p \ \eps^{-2} n^{p/2}  \log n \rfloor
$$
then for any $t > \eps$,
with probabiblity greater than $1 - C \exp (- \left( t / C_p' \eps
\right)^{1/p})$, for any $y \in \R^n$,
$$
(1 - t) \E | \lan X, y \ran |^p \le
 \frac{1}{m} \sum_{j=1}^m |\lan X_j, y \ran|^p \le (1 + t) \E | \lan X, y \ran |^p.
 $$
The constants $C_p$ and $C_p'$ are positive real numbers depending only on $p$.
\end{theorem}
 Theorem \ref{logconcaveintro} provides an almost isometric
approximation of the exact moments, instead of the isomorphic
estimates of \cite{GiaMil}, and achieves it with fewer sample
vectors.
 In
the  case $p=2$, it also improves the estimate of \cite{Rud3}, and
extends to the general setting the estimate obtained by
Giannopoulos, Hartzoulaki and Tsolomitis \cite{GHT} for a random
vector  uniformly distributed in a 1-unconditional isotropic
convex body.

The rest of the paper is organized as follows. In Section
\ref{generalrandomvector} we formulate and prove the main results
for abstract random vectors. The key step of the proof of Theorem
\ref{thrandomintro} is the estimate of  the Gaussian random
process
$$
Z_y = \sum_{j=1}^m g_j | \langle X_j,y \rangle |^p,
$$
where $g_j$ are independent standard Gaussian random variables
${\cal N}(0,1)$. To obtain such estimate we construct  an
appropriate majorizing measure and apply the Majorizing measure
theorem of Talagrand \cite{Tal:genericchaining}. In Sections
\ref{subspace} and \ref{seclogconcave}, we provide applications of
Theorem \ref{thrandomintro}. Since we require only the existence
of high order moments of the norm of $X$ we can apply Theorem
\ref{thrandomintro} to the measures  supported by the contact
points of a convex body, like in \cite{Rud1}, \cite{Rud2}, as well
as to finding a short Lewis-type decomposition, as described in
Section \ref{subspace}.
 In Section \ref{seclogconcave},
we study in detail the case of $\log$-concave random vectors $X$.
In the last part of this paper, we extend the results obtained in
\cite{GiaMil} for a uniform distribution on a discrete cube to a
general random vector $X$, which satisfies a $\psi_2$ estimate for
the scalar products $\pr{X}{y}, \ y \in \R^n$.

%
%
%
%
%
\section{Maximal deviation of the empirical $p$-moment}\label{generalrandomvector}
\subsection{Statement of the results}
Let $K \subset \R^n$ be a convex symmetric body. The modulus of
convexity of $K$ is defined for any $\eps \in (0,2)$ by
$$
  \delta_K (\eps) = \inf \left \{1 - \norm{\frac{x+y}{2}}_K, \|x\|_K = 1,
  \|y\|_K = 1, \|x-y\|_K > \eps \right \}.
$$
We say that $K$ has modulus of convexity of power type $q \ge 2$
if  $\delta_K(\eps) \ge c \eps^q$ for every $\eps \in (0,2)$. It
is  known (see e.g., \cite{Pisier}, Proposition 2.4 or
\cite{Figiel}) that this property is equivalent to the fact that
the inequality
$$
\norm{ \frac{x+y}{2}}_K^q + \la^{-q} \norm{ \frac{x-y}{2}}_K^q \le
\frac{1}{2} (\|x\|_K^q + \|y\|_K^q).
$$
holds for all $x,y \in \R^n$. Here $\la>0$ is a constant depending
only on $c$ and $q$. Referring to this inequality below, we shall
say that $K$ has modulus of convexity of power type $q$ with
constant $\la$.

Our main result is the following theorem, which implies Theorem
\ref{thrandomintro} from the Introduction.

\begin{theorem}\label{thrandom}
Let $K \subset \R^n$ be a symmetric convex body of radius $D$.
 Assume that $K$ has modulus of
convexity of power type $q$ with constant $\la$  for some $q \ge
2$, and let $q^*$ be the conjugate of $q$.
\\
Let $X$ be a random vector in $\R^n$ and let $X_1, \ldots, X_m$ be
independent copies of $X$. For   $p \ge q$ set
$$
A = C^p \, \la^p  \frac{(\log m)^{1/q^*}}{\sqrt m} \dis (D
\gr_{p,m}(X))^{p/2} \hbox{ {\it and} } B= \sup_{y \in K} \E|\lan
X, y \ran|^p.
$$
 Then
$$
 \E
 \sup_{y \in K} \bigg| \frac{1}{m} \sum_{j=1}^m |\lan X_j, y \ran|^p
 - \quad \E |\lan X, y \ran|^p \bigg|
\le
A^2 + A \sqrt{B}.
$$
\end{theorem}
The assumption of Theorem \ref{thrandomintro} reads $A^2 \le \de^2
\cdot B$, hence $A^2+A \sqrt{B} \le 2 \de B$. Thus, Theorem
\ref{thrandomintro} follows immediately from Theorem
\ref{thrandom}.

\noindent
{\bf Remark.} In fact we shall prove a slightly better inequality.
Define
$$
\gr_{p,m}'(X, K) = \bigg( \E \max_{1 \le j \le m} |X_j|_2^2
\
\max_{1 \le j \le m} \|X_j\|_{K^o}^{p-2}
\bigg)^{1/p},
$$
then  Theorem \ref{thrandom} holds, if the quantity $(D
\gr_{p,m}(X))^{p/2} $ is replaced by $D \gr_{p,m}'(X,K)^{p/2}$.
Since $K \subset D B_2^n$, it is clear that
$$
\gr_{p,m}'(X, K)^{p/2}
\le
D^{p/2 -1}
\gr_{p,m}(X)^{p/2}.
$$

The proof of this Theorem  is based on the following
lemma.
%
%
%
%
 \begin{lemma}\label{majorizing}
Let $K \subset \R^n$ be a symmetric convex body of radius $D$.
Assume that $K$ has modulus of convexity of power type $q$ with
constant $\la$ for some $q \ge 2, \la>0$. Let $q^*$ be the
conjugate of $q$.
\\
Then for every $p \ge q$, and every
 deterministic vectors $X_{1}, \ldots, X_{m}$ in $\reel^n,$
$$
\begin{array}{c}
\dis \E \sup_{y \in K}
 \bigg|
 \sum_{j=1}^m \eps_{j} |\lan X_j, y \ran|^p
 \bigg|
\le
\\
\dis C^p \, \la^p \ (\log m)^{1/q^*}      \, D \,
 \max_{1 \le j \le m} |X_{j}|_{2} \,
 \sup_{y \in K} \left(\sum_{j=1}^m |\lan X_j, y
 \ran|^{2(p-1)} \right)^{1/2}
\end{array}
 $$
 where expectation is taken over the Bernoulli random
variables $(\eps_{j})_{1 \le j \le m}$.
 \end{lemma}
 %
%
%
%
 The proof of the Lemma uses a specific construction of a
 majorizing measure. It will be presented in  part \ref{construction}.
 \\
{\bf Proof of Theorem \ref{thrandom}.} The proof is based on a
standard symmetrization argument. We denote by $X'_{1}, \ldots,
X'_{m}$ independent copies of $X_1, \ldots, X_m$. Let
$(\eps_{j})_{j=1}^m$ be independent symmetric Bernoulli random
variables, which are independent of all others.  Then the
expectation of
$$
V_p(K) = \sup_{y \in K} \bigg| \frac{1}{m} \sum_{j=1}^m |\lan X_j, y \ran|^p
 - \quad \E |\lan X, y \ran|^p \bigg|
$$
can be estimated as follows:
$$
\begin{array}{rcl}
m \ \E V_p(K) & =  &
\dis
\E \sup_{y \in K} \bigg|
\sum_{j=1}^m |\lan X_j, y \ran|^p
\ - \ m \, \E |\lan X, y \ran|^p
\bigg|
 \\
& = &
\dis \E \sup_{y \in K}
\bigg|
 \sum_{j=1}^m
\big(|\lan X_j, y \ran|^p - \E |\lan X'_j, y \ran|^p \big)
\bigg|
\\
& \le &
\dis
\E_X \E_{X'} \sup_{y \in K}
\bigg|
\sum_{j=1}^m
\big( |\lan X_j, y \ran|^p - |\lan X_j', y \ran|^p \big)
\bigg|
\\
& = &
\dis
\E_X \E_{X'} \E_{\eps}
\sup_{y \in K}
\bigg|
\sum_{j=1}^m
\eps_j \big( |\lan X_j, y \ran|^p - |\lan X_j', y \ran|^p \big)
\bigg|
\\
& \le &
\dis
2 \E_X \E_{\eps} \sup_{y \in K}
\bigg|
\sum_{j=1}^m \eps_j  |\lan X_j, y \ran|^p
\bigg|.
\end{array}
$$
Therefore,  Lemma \ref{majorizing} implies
$$
\E V_p(K) \le C^p \, \la^p D  \frac{(\log m)^{1/q^*}}{\sqrt m}
\E_{X} \max_{1 \le j \le m} |X_{j}|_{2} \sup_{y \in K} \bigg(
\frac{1}{m}\sum_{j=1}^m |\lan X_j, y \ran|^{2(p-1)} \bigg)^{1/2}.
$$
Since $p \ge 2$,
it is easy to see that
$$
\begin{array}{rl}
\dis
\E_{X}  \max_{1 \le j \le m} |X_{j}|_{2}
\sup_{y \in K}
\bigg( \frac{1}{m}\sum_{j=1}^m |\lan X_j, y \ran|^{2(p-1)} \bigg)^{1/2}
& \le
\\
\dis
\E_{X}  \max_{1 \le j \le m} |X_{j}|_{2}
\max_{1 \le j \le m} \|X_{j}\|_{K^o}^{p/2 -1}
\sup_{y \in K}
\bigg( \frac{1}{m}\sum_{j=1}^m |\lan X_j, y \ran|^{p} \bigg)^{1/2}
& \le
\\
\dis
\gr_{p,m}'(X, K)^{p/2}
\bigg( \E_{X} \sup_{y \in K}
\frac{1}{m}\sum_{j=1}^m |\lan X_j, y \ran|^{p} \bigg)^{1/2}
& \le
\\
\dis
\gr_{p,m}'(X, K)^{p/2}
\bigg( \E V_p(K) + \sup_{y \in K} \E |\lan X, y \ran|^p  \bigg)^{1/2}.
\end{array}
$$
We get that
$\E V_p(K) \le A' ( \E V_p(K) + B)^{1/2}$
where
$$
A' = C^p\, \la^p  D \frac{(\log m)^{1/q^*}}{\sqrt m} \gr_{p,m}'(X,
K)^{p/2} \hbox{ and } B = \sup_{y \in K} \E |\lan X, y \ran|^p
$$
which proves the announced result.
$\hfill \Box$

We present now a  deviation inequality for the positive random
variable $V_p(K)$ under the assumption that $|X|_2$ satisfies some
$\psi_{\alpha}$ estimate. Mendelson and Pajor \cite{MenPa} studied
the same deviation inequality in the case $p=2$ and $K = B_2^n$
using a symmetrization argument. Our approach is based on a
concentration result of Talagrand (Theorem 6.21 in \cite{LT}).
\begin{theorem}\label{concentration}
With the same notation as in Theorem \ref{thrandom}, let $V_p(K)$
be the random variable
$$
V_p(K) = \sup_{y \in K} \bigg| \frac{1}{m} \sum_{j=1}^m |\lan X_j, y \ran|^p
 - \quad \E |\lan X, y \ran|^p \bigg|.
$$
Assume that  $\| |X|_2 \|_{\psi_{\alpha}} < \infty$ for some $0<
\alpha \le p$. Then there exists a positive constant
$c_{\alpha,p}$ depending only on $\alpha$ and $p$ such that
$$
\forall t > 0, \ \Prob (V_p(K) \ge t) \le 2 \exp\left( - ( t /
Q)^{\alpha / p} \right)
$$
where
$$
Q = c_{\alpha,p} \left( \E V_p(K) + \frac{(\log m)^{p/\alpha}}{m} D^p \| |X|_2 \|_{\psi_{\alpha}}^p \right).
$$
\end{theorem}
{\bf Remark.}
Observe that in the typical case,  $Q$ is of the order $\E
V_p(K)$ for which we may use Theorem \ref{thrandom}. By Lemma \ref{psia} (see below),
$$
\gr_{p,m} (X) \le C (p \, \log m)^{1/\al} \| |X|_2 \|_{\psi_{\alpha}}
$$
therefore, using Theorem \ref{thrandom},
$$
Q \le C_{\al,p} (2 A_1^2 + A_1 \sqrt B)
$$
where
$$
A_1 = \la^p D^{p/2} \frac{(\log m)^{1/q^* + p/2 \al}}{\sqrt m} \| |X|_2 \|_{\psi_{\alpha}}
\hbox{ and }
B = \sup_{y \in K} \E |\lan X,y\ran |^p.
$$
For the proof of this theorem, we need an elementary lemma.
\begin{lemma}\label{psia}
  Let $\delta>0$ and let $Z_1, \ldots, Z_m$ be independent copies of
  a random variable $Z$. Then
 $$
 \| \,  \max_{j=1, \ldots, m} |Z_j| \ \|_{\psi_{\delta}}
 \le C \log^{1/\de} m \cdot \| Z \|_{\psi_{\de}}.
 $$
\end{lemma}
{\bf Proof.}
Note that for any random variable $Y$ the inequality
$\norm{Y}_{\psi_\de} \le A$ is equivalent to
\[
 \norm{Y}_r \le C A r^{1/\de}
\]
for all $r >1$. Assume that $r< \log m$. Then
$$
\begin{array}{rcl}
 \dis \| \, \max_{j=1 \ldots m} |Z_j| \ \|_r
  &
  \le &  \dis
   \|  ( \sum_{j=1}^m |Z_j|^{\log m} )^{1/\log m} \|_r
  \le
  \left (
  \E ( \sum_{j=1}^m |Z_j|^{\log m})^{r/\log m}
      \right )^{1/r}
      \\
  &\le & \dis  \left ( \sum_{j=1}^m \E |Z_j|^{\log m}
       \right )^{1/\log m}
  \le C \log^{1/\de} m \cdot \norm{Z}_{\psi_\de}.
\end{array}
$$
If $r> \log m$, then using $\max_{j=1, \ldots, m} a_j \le
(\sum_{j=1}^m a_j^r)^{1/r}$, we get
$$
 \| \, \max_{j=1, \ldots, m} |Z_j| \ \|_r
 \le ( \sum_{j=1}^m \E |Z_j|^r )^{1/r}
 \le m^{1/r} \norm{Z}_r
 \le C r^{1/\de} \cdot \norm{Z}_{\psi_\de}.
$$
These two inequalities imply the Lemma.  $\hfill \Box$

 {\bf Proof of Theorem
\ref{concentration}.} To any vector $x \in \R^n$ we associate the
function $f_x$ defined on $K$ by
$$
\begin{array}{rcl}
f_x  :  K & \to & \R
\\
y & \mapsto &  \frac{1}{m} \left( |\lan x, y \ran|^p - \E |\lan X,y\ran |^p \right).
\end{array}
$$
Let $f_X$ be the random vector of $L_{\infty}(K)$ associated to $X$.
Now we apply Theorem 6.21 of Ledoux--Talagrand \cite{LT} to $\norm{
\sum_{j=1}^m f_{X_j}}$ where the $X_j$'s are independent copies of $X$.
By definition,
$$
\| \sum_{j=1}^m f_{X_j} \|_{L_{\infty}(K)} = V_p(K),
$$
and
$$
\| f_{X} \|_{L_{\infty}(K)}
\le
 \frac{1}{m} \left( \sup_{y\in K} |\lan X,y \ran|^p +
\sup_{y\in K} \E|\lan X,y \ran|^p \right)
\le
\frac{D^p}{m} \left(  |X|_2^p + \E |X|_2^p \right).
$$
Theorem 6.21 of Ledoux-Talagrand \cite{LT} states that if
$\alpha / p \le 1$, there exists a constant $c_{\alpha,p}$
depending only on $\alpha/p$ such that
$$
\| V_p(K) \|_{\psi_{\alpha/p}}
\le
c_{\alpha,p} \left(
\E V_p(K) + \| \max_{1 \le j \le m} \|f_{X_j}\|_{L_{\infty}(K)} \, \|_{\psi_{\alpha/p}}
\right).
$$
Moreover,
\begin{align*}
 \| \max_{1 \le j \le m} \|f_{X_j}\|_{L_{\infty}(K)} \,
  \|_{\psi_{\alpha/p}}
  &\le  \frac{2 D^p}{m} \ \| \max_{1 \le j \le m}
  |X_j|_2^p \, \|_{\psi_{\alpha/p}}
  \\
  &=
  \frac{2 D^p}{m} \ \| \max_{1 \le j \le m} |X_j|_2 \,
  \|_{\psi_{\alpha}}^p.
\end{align*}
Lemma \ref{psia} implies
\begin{equation}\label{replace}
\| \max_{1 \le j \le m} |X_j|_2 \,  \|_{\psi_{\alpha}}
\le
C (\log m)^{1/\alpha} \| \ |X|_2 \|_{\psi_{\alpha}}.
\end{equation}
This proves that
$$
\| V_p(K) \|_{\psi_{\alpha/p}}
\le
c_{\alpha,p} \left(
\E V_p(K) +  \frac{(\log m)^{p/\alpha}}{m} \, D^p \,  \| \ |X|_2 \|_{\psi_{\alpha}}^p \right).
$$
The deviation inequality follows from the Chebychev inequality.
$\hfill \Box$
\subsection{Construction of majorizing measures}\label{construction}
%
%
%
%
Let us recall the assumptions of Lemma \ref{majorizing}. The
ambient space is $\reel^n$ equipped with a Euclidean structure and
we denote by $|\cdot|_2$ the norm associated. The symmetric convex
body $K$ has a modulus of convexity of power type $q \ge 2$ with a
constant $\la$,  which means that
\begin{eqnarray}\label{Clark}
\forall x,y \in \reel^n, \norm{\frac{x+y}{2}}_K^q + \la^{-q}
\norm{ \frac{x-y}{2}}_K^q \le \frac{1}{2} (\|x\|_K^q + \|y\|_K^q).
\end{eqnarray}
and satisfies also the inclusion $K \subset D B_2^n$, which
means that
$$
\forall x \in \reel^n,   |x|_2 \le D \|x\|_K.
$$
Let $p \ge q \ge 2$, and $X_1, \ldots, X_m$ be $m$ fixed vectors in
$\reel^n$. We define the random process $V_y$  for all $y \in
\reel^n$ by
$$
V_y = \sum_{j=1}^m \eps_j | \langle X_j,y \rangle |^p,
$$
where $\eps_j$ are independent symmetric Bernoulli random
variables. It is well known that this process satisfies a
sub-Gaussian tail estimate: $\forall y, \overline y \in \reel^n$,
$\forall t > 0$,
$$
P(|V_y - V_{\overline y}| \ge t)
\le
2
\exp\left(- \frac{c t^2}{\tilde{d}^2(y,\overline y)} \right)
$$
where
$$
\tilde{d}^2(y,\overline y) =
\sum_{j=1}^m  \left(
| \langle X_j,y \rangle |^p
- | \langle X_j, \overline y \rangle |^p
\right)^2.
$$
Instead of working with this function which is not a metric, it
will be preferable to consider the following quasi-metric
$$
d^2(y,\overline y) =
\sum_{j=1}^m  | \langle X_j,y-\overline y \rangle |^{2}
\left(
| \langle X_j,y \rangle |^{2(p-1)}
+ | \langle X_j, \overline y \rangle |^{2(p-1)}
\right).
$$
The following propositions state inequalities
that we will need to prove Lemma \ref{majorizing}.
Proposition \ref{prop:geometry}  gives some information
concerning the
geometry of the balls associated to the metric $d$ and
Proposition \ref{prop:d,sup,eucl}
explains relation between metric $d$,  new
Euclidean norm and the following
norm defined by
$$
\|x\|_{\infty} = \max_{1 \le j \le m} |\lan X_j,x\ran|.
$$
We denote by $\B_{\rho} (x)$ the ball of center $x$ with radius
$\rho$ for the quasi-metric $d$.
%
%
%
%
%
%
%
%
\begin{proposition}\label{prop:geometry}
For all $y, \overline y \in K$
\begin{eqnarray}
\label{ineg1}
\tilde{d}(y, \overline y) \le p \, d(y,\overline y),
\\
\label{ineg2}
d(y,\overline y)
\le  \sqrt 2 \|y-\overline y\|_{\infty} \
\sup_{y \in K}
\left( \sum_{j=1}^m
| \langle X_j,y \rangle |^{2(p-1)}
\right)^{1/2},
\\
\label{ineg3}
\|y - \overline y\|_{\infty} \le D \max_{1 \le j \le m} |X_j|_2
\ \|y-\overline y\|_K.
\end{eqnarray}
Moreover, the quasi-metric $d$ satisfies the generalized triangle
inequality, and for any point $x$, the ball $\B_{\rho} (x)$ is a
convex set: for all $u_1, \ldots, u_N \in \reel^n$ and all $x, y,
z \in \reel^n$,
\begin{eqnarray}\label{ineg:triangle}
d(u_1, u_N) \le 2p \sum_{i=1}^{N-1} d(u_i, u_{i+1})
\hbox{ and }
d^2(x, \frac{y+z}{2}) \le \frac{1}{2} \big(d^2(x,y) + d^2(x,z) \big).
\end{eqnarray}
\end{proposition}
%
%
To prove it, we will need the following
basic inequalities on real numbers.
\begin{lemma}\label{lem:dist}
For every $x, y \in \reel^+$ and $p\ge 2$, we have
\begin{eqnarray}\label{ineg:classic}
|x^p - y^p| \le p | x-y | \sqrt{x^{2p-2} + y^{2p-2}}
\end{eqnarray}
Moreover, if  $f(s,t) = |s-t| \sqrt{|s|^{2p-2} + |t|^{2p-2}}$ then
for all $r_1, \ldots, r_N \in \reel$
$$
f(r_1, r_N) \le 2p \sum_{i=1}^{N-1} f(r_i, r_{i+1})
$$
and for all $r,$ $s$, $t \in \reel$,
$$
f(r,(s+t)/2)^2 \le \big(f(r,s)^2 + f(r,t)^2\big)/2
$$
\end{lemma}
{\bf Proof.} The first inequality is straightforward. To prove the
second one, consider two cases. When $r_1 r_N \ge 0$, since
$$
  |r_1-r_N| \sqrt{|r_1|^{2p-2} + |r_N|^{2p-2}}
  \le \sqrt 2 \Big ||r_1|^p - |r_N|^p \Big |,
$$
the conclusion follows from the triangle inequality and inequality
(\ref{ineg:classic}). When $r_1 r_N \le 0$, we can assume without
loss of generality that $r_1 \ge 0$ and $r_N \le 0$. Then
\begin{align*}
  f(r_1, r_N)
  &= (r_1+|r_N|) \sqrt{r_1^{2p-2} + |r_N|^{2p-2}}
  \le
  (r_1+|r_N|)(r_1^{p-1} + |r_N|^{p-1}) \\
  &\le 2(r_1^p + |r_N|^{p}).
\end{align*}
 Let $m<N$ be a number such that $r_m \ge 0$ and $r_{m+1} \le 0$.
Then
\[
  r_1^p + |r_N|^{p}
  \le \sum_{i=1}^{m-1} \Big | |r_i|^p-|r_{i+1}|^p \Big |
      +r_{m}^p + |r_{m+1}|^{p}
      + \sum_{i=m+1}^{N-1} \Big | |r_i|^p-|r_{i+1}|^p \Big |.
\]
Combining the previous inequalities with  (\ref{ineg:classic}), we
get
\begin{align*}
  f(r_1, r_N)
  & \le  2\sum_{i=1}^{m-1} \Big | |r_i|^p-|r_{i+1}|^p \Big |
      +2(r_{m}^p + |r_{m+1}|^{p})
      + 2\sum_{i=m+1}^{N-1} \Big | |r_i|^p-|r_{i+1}|^p \Big |
\\
  & \le 2p \sum_{i=1}^{m-1} f(r_i, r_{i+1})
    + 2(r_{m} - r_{m+1})(r_{m}^{p-1} + |r_{m+1}|^{p-1}) \\
  & + 2p \sum_{i=m+1}^{N-1} f(r_i, r_{i+1})
\\
  & \le
   2p \sum_{i=1}^{m-1} f(r_i, r_{i+1})
     + 2 \sqrt 2 f(r_{m}, r_{m+1})
     + 2p \sum_{i=m+1}^{N-1} f(r_i, r_{i+1})
\end{align*}
which proves the announced result. The last inequality follows
from the fact that for $p \ge 2,$ the function $v \mapsto
(1-v)^2(1+v^{2p-2})$ is convex on $\reel$, which can be checked by
computing the second derivative. $\hfill \Box$

\medskip
\noindent {\bf Proof of Proposition \ref{prop:geometry}.}
Inequalities (\ref{ineg1}) and (\ref{ineg:triangle}) clearly
follow from the three inequalities proved in Lemma \ref{lem:dist}.
Inequalities (\ref{ineg2}) and (\ref{ineg3}) follow from simple
observations about $d$ and the fact that $K \subset D B_2^n$.
$\hfill \Box$
%
\begin{proposition}\label{prop:d,sup,eucl}
Let
$
\dis
M = \sup_{y \in K}\sum_{j=1}^m
| \langle X_j,y \rangle |^{2(p-1)}.
$
For a fixed $u \in K$,  we define the Euclidean norm $|\cdot|_{{\cal E}_u}$
associated to $u$ by
$$
|z|_{{\cal E}_u}^2 =
\sum_{\ell = 1}^m
| \langle X_{\ell}, z  \rangle |^2
| \langle X_{\ell}, u \rangle |^{2(p-1)},
\ \forall z \in \reel^n.
$$
Then the following inequality holds for all $z$,
 $\overline z \in \reel^n$:
$$
d^2(z, \overline z) \le
2 \cdot 4^{p-1} \big( |z-\overline z|_{{\cal E}_u}^2 +
 M \ \|z-\overline z\|_{\infty}^2
(\|z-u\|_K^{2p-2}
+\|\overline z -u \|_K^{2p-2} )\big).
$$
\end{proposition}
%
%
%
%
%
{\bf Proof.} By homogeneity of the statement, we can assume that
$$
M = \sup_{y \in K}\sum_{j=1}^m
| \langle X_j,y \rangle |^{2(p-1)}=1.
$$
For any $z \in \reel^n$, let $L_z = \{\ell \in \{1, \ldots, m\}\,
\big | \ |\lan X_{\ell}, z \ran | \ge 2 |\lan X_{\ell}, u \ran |
\}$. Then by convexity of the function $t \mapsto t^{2p-2}$, we
have
$$
\begin{array}{rcl}
\dis \sum_{\ell \in L_z}
| \langle X_{\ell}, z \rangle |^{2(p-1)}
& \le &
\dis 2^{2p-3} \sum_{\ell \in L_z}
| \langle X_{\ell}, z - u  \rangle |^{2(p-1)}
+
2^{2p-3} \sum_{\ell \in L_z}
| \langle X_{\ell}, u  \rangle |^{2(p-1)}
\\
& \le & \dis 2^{2p-3} \sum_{\ell \in L_z} | \langle X_{\ell}, z -
u  \rangle |^{2(p-1)} + \frac{1}{2} \sum_{\ell \in L_z} | \langle
X_{\ell}, z \rangle |^{2(p-1)},
\end{array}
$$
which proves (since $M=1$) that for any $z \in \reel^n$,
$$
\sum_{\ell \in L_z}
| \langle X_{\ell}, z \rangle |^{2(p-1)}
\le 4^{p-1} \|z - u\|_K^{2p-2}.
$$
Hence, for any $z, \overline z \in \reel^n$,
\begin{align*}
   \sum_{\ell \in L_z} | \langle X_{\ell}, z - \overline z
    \rangle |^2 | \langle X_{\ell}, z \rangle |^{2(p-1)}
  & \le  \dis\|z - \overline z\|_{\infty}^2
     \sum_{\ell \in L_z}| \langle X_{\ell}, z \rangle |^{2(p-1)}
    \\
  & \le
     4^{p-1} \|z - u\|_K^{2p-2}  \, \|z- \overline z
     \|_{\infty}^2.
\end{align*}
For any $l \notin L_z$ we have $|\pr{X_l}{z}| \le 2
|\pr{X_l}{u}|$, so
\[
    \sum_{\ell \notin L_z} | \langle X_{\ell}, z -
   \overline z \rangle |^2 | \langle X_{\ell}, z \rangle |^{2(p-1)}
    \le 4^{p-1} \sum_{\ell = 1}^m | \langle X_{\ell}, z - \overline z
   \rangle |^2 | \langle X_{\ell}, u \rangle |^{2(p-1)}.
\]
The same inequalities hold if we exchange the roles of $z$ and
$\overline z$. To compute  $d^2(z_i, z_j)$, we split the  sum in
four parts and apply the inequalities above:
$$
\begin{array}{l}
d^2(z, \overline z)  =
\dis
\sum_{\ell = 1}^m
| \langle X_{\ell}, z - \overline z  \rangle |^2
| \langle X_{\ell}, z \rangle |^{2(p-1)}
+
| \langle X_{\ell}, z - \overline z  \rangle |^2
| \langle X_{\ell}, \overline z \rangle |^{2(p-1)}
\\
 =
\dis
\sum_{\ell \in L_z}
| \langle X_{\ell}, z - \overline z  \rangle |^2
| \langle X_{\ell}, z \rangle |^{2(p-1)}
+ \sum_{\ell \notin L_z}
| \langle X_{\ell}, z - \overline z  \rangle |^2
| \langle X_{\ell}, z \rangle |^{2(p-1)}
\\
 + \dis
\sum_{\ell \in L_{\overline z}}
| \langle X_{\ell}, z - \overline z  \rangle |^2
| \langle X_{\ell},  \overline z \rangle |^{2(p-1)}
+ \sum_{\ell \notin L_{\overline z}}
| \langle X_{\ell}, z - \overline z  \rangle |^2
| \langle X_{\ell}, \overline z \rangle |^{2(p-1)}
\\
\le \dis 2 \cdot 4^{p-1} \left(|z - \overline z |_{{\cal E}_u}^2 +
\|z- \overline z \|_{\infty}^2 (\|z - u\|_K^{2p-2} + \|z -
u\|_K^{2p-2}) \right). \hfill \Box
\end{array}
$$
%
%
{\bf Proof of Lemma \ref{majorizing}.} By inequality
(\ref{ineg1}),  we may treat $V_y$ as a sub-Gaussian process with
the quasi-metric $p \cdot d$. By homogeneity of the statement, we
can assume that
$$
\sup_{y \in K} \sum_{j=1}^m | \langle X_j,y \rangle |^{2(p-1)} =
1.
$$
Denote  $Q = \max_{1\le j \le m} |X_j|_2$. We want to show that
\begin{eqnarray}\label{sup}
\E \sup_{y \in K} |V_y| \le C^p \, \la^p \ Q \ (\log m)^{1/q^*} \
D.
\end{eqnarray}
By Proposition \ref{prop:geometry}, the diameter of the set $K$
with respect to the metric $d$ is bounded by $2\sqrt 2 Q D$. Let
$r$ be a fixed number chosen such that $r = c p^2$ for a large
universal constant $c$ and $k_0$ be the largest integer such that
$r^{-k_0} \ge 2\sqrt 2 QD$.

The proof of inequality $(\ref{sup})$ is based on the majorizing
measure theory of Talagrand \cite{Tal:genericchaining}. The
following theorem is a combination of Proposition 2.3, Theorem 4.1
and Proposition 4.5 of \cite{Tal:genericchaining}. Note that
assuming that $r \ge 2$, one can set $K(r)=C$ in Proposition 2.3,
and $K(2,1,r)=C$ in Proposition 4.5.

\begin{thm}[\cite{Tal:genericchaining}]
{\it Let $r \ge 2$. Let  $\phi_k : K \to \reel^+$ for $k \ge k_0$
be a family of maps satisfying the following assumption: there
exists $A>0$ such that for any point $x \in K$, for any $k \ge
k_0$ and any $N \in \entier$
\\
$
(H)
\left\{
\begin{array}{l}
for \ any \ points \ x_1, \ldots, x_N \in \B_{r^{-k}}(x)
\ with \ d(x_i, x_j) \ge r^{-k-1}, i \ne j
\\
we \ have \ \dis \max_{i=1, \ldots, N} \phi_{k+2}(x_i)
\ge \phi_k(x) + \, \frac{1}{A} \ r^{-k} \sqrt{\log N}.
\end{array}
\right.
$
\\
Then for any fixed $y_0 \in K$,
$$
  \E \sup_{y \in K} |V_y - V_{y_0}|
  \le c \ A \cdot \sup_{k\ge k_0, x\in K} \phi_k(x).
$$
}
\end{thm}
%
%

\bigskip
To obtain  the conclusion of Lemma \ref{majorizing}, set $y_0 =
0$.
\\
To complete the proof , we have to define  the functionals $\phi_k
: K \to \reel^+$. Let $k_1$ be the smallest integer such that
$r^{-k_1} \le QD / \sqrt n$. For $k \ge k_1 +1$, set

$$
\phi_k(x) =
1 + \frac{1}{2 \log r} + \frac{\sqrt n}{Q \, D \, (\log m)^{1/q^*}}
\sum_{l=k_1}^k r^{-l}
\sqrt{\log(1 + 4 Q D  r^{l})}.
$$
Note that in this range of $k$ the functionals $\phi_k$ do not
depend on $x$. We shall show that with this choice of $\phi_k$,
the condition $(H)$ follows from the classical volumetric estimate
of the covering numbers.

For $k_0 \le k \le k_1$, the functionals $\phi_k$ are defined by
$$
\phi_k(x) =
\min\{ \|y\|_K^q, y \in \B_{4pr^{-k}}(x) \} +
\frac{k-k_0}{\log m}.
$$
Since $q \ge 2$ then $1\le q^* \le 2$ and
$(\log m)^{1/q^*} \ge \sqrt{\log m}$.
It is easy to see using
definitions of $k_0$ and $k_1$ that
$$
\sup_{x \in K, k\ge k_0} \phi_k(x) \le c.
$$
We shall prove that our functionals satisfy condition $(H)$ for
$$
A = (C \la )^p \ Q \ D \ (\log m)^{1/q^*}
$$
where $C$ is a large numerical constant.
That will conclude
the proof of Lemma \ref{majorizing} with a new constant $C$.
$\hfill \Box$
\\
{\bf Proof of condition $(H)$.}
 Let $N \in \entier$, $x \in K$,
$x_1, \ldots, x_N \in \B_{r^{-k}}(x)$ with $d(x_i, x_j) \ge
r^{-k-1}$. We have to prove that
$$
\max_{i=1, \ldots, N} \phi_{k+2} (x_i) - \phi_k(x)
 \ge \frac{r^{-k} \sqrt{\log N}}{ (C \lambda)^p \ Q \ D \ (\log
m)^{1/q^*}}.
$$

For $k \ge k_1-1$, we always have
$$
\phi_{k+2}(x_i) - \phi_k(x)
\ge
\frac{\sqrt {n \log(1 + 4 Q D r^{k+2})}}
{Q D (\log m)^{1/q^*}} \ r^{-k-2}.
$$
Since the points $x_1, \ldots, x_N$ are well separated in the
metric $d$, they are also well separated in the norm
$\|\cdot\|_K$. Indeed, by $(\ref{ineg2})$ and $(\ref{ineg3})$, we
have
\[
 \|x_i - x_j\|_K \ge r^{-k-1} / Q D \sqrt 2.
\]
  By the classical
 volumetric estimate, the maximal cardinality of a $t$-net in a
 convex symmetric body $K \subset \R^n$ with respect to $\|\cdot
 \|_K$ does not exceed
$
    \left ( 1+ 2/t \right )^n.
$
 Therefore,
$$
\sqrt{\log N}  \le \sqrt{n \log(1 + 2 \sqrt{2} QD r^{k+1})},
$$
which proves the desired inequality.

The case $k_0 \le k \le k_1-2$ is much more difficult. Our proof
uses estimates of the covering numbers, in particular, the dual
Sudakov inequality \cite{PajTom}. Recall that the covering number
$N(W,\|\cdot \|_X,t)$ is the minimal cardinality of
$\|\cdot\|_X$-balls of radius $t$ needed to cover the $W$.

For $j = 1, \ldots, N$  denote by $z_j \in K$  the points which
satisfy $ \|z_j\|_K^q = \min \{ \|y\|_K^q, y \in
\B_{4pr^{-k-2}}(x_j) \} $. Denote by $u \in K$ a point such that $
\|u\|_K^q = \min \{ \|y\|_K^q, y \in \B_{4pr^{-k}}(x) \}. $ Set
$$
\theta = \max_{j} \|z_j\|_K^q - \|u\|_K^q.
$$
Then we  have $ \max_{j} \phi_{k+2} (x_j) - \phi_k(x) = \theta +
\frac{2}{\log m} $. We shall prove that
\begin{eqnarray}\label{wanted}
\theta + \frac{2}{\log m}
\ge
r^{-k} \sqrt{\log N} / A.
\end{eqnarray}
Since $d(x_i, x_j) \ge r^{-k-1}$, $z_l \in \B_{4pr^{-k-2}}(x_l)$,
and $d$ satisfies a generalized triangle inequality,
the points $(z_j)_{1\le j \le N}$ remain well separated. Indeed,
$$
r^{-k-1}
\le
d(x_i, x_j) \le 2p (d(x_i, z_i) + d(z_i, z_j) + d(z_j, x_j))
\le
2p d(z_i, z_j) + 16 p^2 r^{-k-2}
$$
and since $r=cp^2$, we have
$$
d(z_i, z_j) \ge r^{-k-1}/cp
$$
for all $i \ne j$. Recall that $r=cp^2$. Using again the
generalized triangle inequality, we get that
$$
d(x,z_j) \le 2p(d(x, x_j) + d(x_j, z_j)) \le 2p (r^{-k} + 4p
r^{-k-2}) \le 4p r^{-k}.
$$
It means that $z_j \in \B_{4pr^{-k}}(x)$, $u \in
\B_{4pr^{-k}}(x)$, and the convexity of the balls for the
quasi-metric $d$ proved in Proposition \ref{prop:geometry} implies
 $(u+z_j)/2 \in \B_{4pr^{-k}}(x)$. Since $K$ has modulus of
convexity of power type $q$, inequality $(\ref{Clark})$ holds. By
the definition of $u$, we get that for all $j=1, \ldots, N$
$$
\la^{-q} \norm{\frac{z_j - u}{2}}_K^q \le \frac{1}{2} \left(
\|z_j\|_K^q + \|u\|_K^q \right) - \norm{\frac{z_j + u}{2}}_K^q \le
\frac{\|z_j\|_K^q - \|u\|_K^q}{2} \le \frac{\theta}{2}.
$$
This proves that $\forall j=1, \ldots, N, \|z_j - u \|_K \le 2 \la
\theta^{1/q}. $ Let $\delta>0$. Consider the set
$$
U = u + 2 \la \theta^{1/q} K
$$
which contains all the $z_j$'s and let $S$ be the maximal number
of points in $U$ that are $2 \delta$ separated in
$\|\cdot\|_{\infty}$. Then $U$ is covered by $S$  subsets of
diameter smaller than $2\delta$ in $\|\cdot\|_{\infty}$ metric,
and so $ S \le N(U, \|\cdot\|_{\infty}, 2 \delta). $ Set
$$
\delta = \tilde{c}^p \la^{1-p} r^{-k} \theta^{1/q - 1}
$$
where the  constant $\tilde{c}$ will be chosen later. Since $U = u
+ 2 \la \theta^{1/q} K$ and $K \subset D B_2^n$,the dual Sudakov
inequality \cite{PajTom} implies
$$
\sqrt{\log S} \le \sqrt{\log N(B_2^n, \|\cdot\|_{\infty}, \delta /
D \la \theta^{1/q})} \le c \ D \ \la \ \theta^{1/q} \ \E
\|G\|_{\infty} \, / \, \delta.
$$
Here $G$ denotes a standard Gaussian vector in $\reel^n$. It is
well known that
$$
\E \|G\|_{\infty} = \E \max_{j=1, \ldots, m} |\lan X_j, G \ran|
\le c \ Q \ \sqrt{\log m}.
$$
We consider now two cases.

First, assume that $S \ge \sqrt N$. Then by previous estimate and
the definition of $\delta$, we get
$$
\sqrt{\log N}
\le c \ Q \ \la \ D \ \sqrt{\log m} \ \theta^{1/q} \, / \, \delta
\le \theta \ \tilde{c}^p
\ r^k \ Q \ \la^p \ D \ \sqrt{\log m}
$$
which easily proves  $(\ref{wanted})$ (since $q^* \le 2$).

The second case is when $S\le \sqrt N$. Since $U$ is covered by
$S$ balls of diameter smaller than $2\delta$ in
$\|\cdot\|_{\infty}$, there exists a subset $J$ of $\{1, \ldots,
N\}$ with $\# J \ge \sqrt N$ such that
$$
\forall i,j \in J, \|z_i - z_j\|_{\infty} \le 2 \delta.
$$
By Proposition \ref{prop:d,sup,eucl} applied to the Euclidean norm
defined  by
$$
|y|_{{\cal E}_u}^2 = \sum_{\ell = 1}^m
| \langle X_{\ell}, y  \rangle |^2
| \langle X_{\ell}, u \rangle |^{2(p-1)},
$$
we get that
$$
d^2(z_i, z_j)
\le
2 \cdot 4^{p-1} \left(|z_i - z_j|_{{\cal E}_u}^2 + 4^p \la^{2p -2} \theta^{(2p-2)/q}
 \delta^2 \right).
$$
Since  $\theta \le 1$ and $q \le p$, the definition of $\delta$
implies
$$
4^{2p} \la^{2p -2} \theta^{(2p-2)/q} \, \delta^2
\le (4\tilde{c})^{2p} r^{-2k} \theta^{2(p/q -1)}
\le (4\tilde{c})^{2p} r^{-2k}.
$$
Recall that $d(z_i, z_j) \ge  r^{-k-1} / cp$ and $r = c p^2$.
Hence,
$$
r^{-2k} / cp^6 \le d(z_i, z_j)^2
 \le 2 \cdot 4^{p-1} |z_i - z_j|_{{\cal E}_u}^2
+ 2  (4\tilde{c})^{2p} r^{-2k}.
$$
Choosing $\tilde{c}$ small enough, we get that for all  $i, j \in
J$,
$$
|z_i - z_j|_{{\cal E}_u} \ge r^{-k-1} c^p.
$$
Since $K \subset D B_2^n$, we have the following  estimate for the
covering numbers:
$$
\begin{array}{rcl}
\# J \le N(U, |\cdot|_{{\cal E}_u}, c^p r^{-k-1})
& = &
N(K, |\cdot|_{{\cal E}_u}, c^p r^{-k-1} / 2 \la \theta^{1/q})
\\
& \le &
N(B_2^n, |\cdot|_{{\cal E}_u}, c^p r^{-k-1}/ 2 \la \theta^{1/q} D).
\end{array}
$$
 Recall that
$G$ denotes a standard Gaussian vector in $\reel^n$. By the dual
Sudakov inequality \cite{PajTom}, we have
$$
\begin{array}{rcl}
\dis
\sqrt{\log N(B_2^n, |\cdot|_{{\cal E}_u},
\frac{c^p r^{-k-1}}{2 \la D \theta^{1/q}})}
& \le &
C^p \ r^{k+1} \ \theta^{1/q} \ \la \ D \
\E |G|_{{\cal E}_u}
\\
& \le & \dis C^p \ r^{k+1} \ \theta^{1/q} \ \la \ D \left( \E
|G|_{{\cal E}_u}^2\right)^{1/2}.
\end{array}
$$
Since for all $y \in K$, $\sum_{j=1}^m | \langle X_j,y \rangle
|^{2(p-1)} \le 1$, we obtain
$$
\E |G|_{{\cal E}_u}^2 = \sum_{\ell = 1}^m |X_{\ell}|_2^2
| \langle X_{\ell}, u \rangle |^{2(p-1)} \le Q^2.
$$
Since $\# J \ge \sqrt N$, we have $ \sqrt{\log N} \le C^p  r^{k+1}
 \la D  Q  \theta^{1/q} $ with a universal constant $C$.
Moreover, by Young's inequality
$$
\theta^{1/q} \le (\log m)^{1/q^*}
\left( \theta / q + 1 / (q^* \log m) \right)
$$
and since $q^* \le 2 \le q$ and $\la \ge 1$, we get
$$
\sqrt{\log N}
\le
(C \la)^p \ r^{k+1} \ D \  Q \
(\log m)^{1/q^*}
\left( \theta + \frac{2}{\log m} \right).
$$
This completes the proof of $(\ref{wanted})$ and the proof of
condition $(H)$ for the functionals $\phi_k$. $\hfill \Box$

\section{Approximate Lewis decomposition}\label{subspace}
%
%
%
%
%
It is well known that if $E$ is an $n$-dimensional subspace of
$L_p$, then $E$ is $(1+\eps)$-isomorphic to an $n$-dimensional
subspace of $\ell_p^N$ with $N$  depending on $n$, $p$ and $\eps$.
Lewis \cite{Lewis} proved that any linear subspace $E$ of
$\ell_p^N$ possesses a special decomposition of the  identity.
More precisely,
 there exists a Euclidean structure on $E$ with
 the scalar product $\langle \cdot, \cdot
\rangle$, vectors $y_{1}, \ldots, y_{N} \in E$ and scalars $c_{1},
\ldots, c_{N} > 0$ such that
$$
\left \{
\begin{array}{l}
\forall i, \pr{y_i}{y_i} = 1 ,
\\
\| x \|_{E} = \dis \bigg( \sum_{i=1}^{N} c_{i} |\langle x, y_{i} \rangle|^p
\bigg)^{1/p} , \forall x \in E,
\\
{\rm Id}_{E} = \dis  \sum_{i=1}^{N} c_{i} y_{i} \otimes y_{i}.
\end{array}
\right.
$$

Denote by $(H, | \cdot |_{H})$ the linear space $E$ equipped with
this Euclidean structure. Recall that $p^*$ denotes the conjugate
of $p$.
 In the following Theorem, we prove that both spaces $E$ and $H$
can be $(1+ \eps)$-embedded in $\ell_p^m$ and  $\ell_2^m$
respectively via the same linear operator $T: \R^N \to \R^m$,
whenever $m$ is of the order of $\eps^{-2} n^{p/2} \log^{2/p^*} (n / \eps^{4/p})$.
This extends a classical result of Bourgain, Lindenstrauss and
Milman \cite{BouLinMil} (and \cite{LT} for a better dependance on
$\eps$) and some results in \cite{Rud1} concerning the number of
contact points of a convex body needed to approximate the identity
decomposition.
\begin{theorem}\label{subLp}
Let $E$ be an $n$-dimensional subspace of $L_p$ for some $p \ge
2$. Then for every $\eps > 0$
 there exists a Euclidean structure $H= (E,\lan \cdot, \cdot \ran)$
on $E$ and $m$ points $x_1, \ldots, x_m$ in $E$ with
$$
m \le \frac{C^p}{\eps^{2}} \ n^{p/2} \ \log^{2/p^*} \left (
\frac{n}{\eps^{4/p}} \right )
 \le \frac{C^p}{\eps^{2}} \ n^{p/2} \log^2 \left (
\frac{n}{\eps^{4/p}} \right )
$$
such that
$ \forall j, | x_{j} |_{H} = 1 $ and for all $y \in E$,
$$
\left \{
\begin{array}{l}
(1 - \eps) \| y \|_{E}
\le
\dis \bigg( \frac{n}{m} \sum_{j=1}^{m}  |\langle y, x_{j} \rangle|^p
\bigg)^{1/p}
\le (1 + \eps) \| y \|_{E}
\\
(1 - \eps) | y |_{H}
\le
\dis \bigg( \frac{n}{m} \sum_{j=1}^{m}  |\langle y, x_{j} \rangle|^2
\bigg)^{1/2}
\le (1 + \eps) | y |_{H}.
\end{array}
\right.
$$
\end{theorem}
{\bf Proof.}
Let $X$ be the random vector taking values $y_{i}$ with probability
$c_{i} / n$. Then  for all $y \in E$,
$$
\E |\langle X, y \rangle|^p =  \| y \|_{E}^p / n
\quad
\hbox{ and }
\quad
\E |\langle X, y \rangle|^2 =   |y|_{H}^{2} / n
\quad
\hbox{ and }
\quad
|X|_{H} = 1.
$$

We will apply Theorem \ref{thrandomintro} twice: first time for
the unit ball of $E$, and then for the unit ball of $H$.

  Since $E$ is a subspace of $L_p$, by Clarkson's
inequality \cite{Clarkson}, $B_E$ has modulus of convexity of
power type $p$ with constant $\la = 1$. From Lewis decomposition,
we get $\| y \|_{E} \le |y|_{H} \le n^{\frac{1}{2} - \frac{1}{p}}
\| y \|_{E}$ which means that for $ D = n^{\frac{1}{2} -
\frac{1}{p}}$,
$$
B_H \subset B_E \subset
D B_H.
$$
Let $X_{1}, \ldots, X_{m}$ be independent copies of $X$, then
$$\sup_{y \in B_E}
 \E|\lan X, y \ran|^p  = 1 / n
\quad
\hbox{ and }
\quad
\gr_{p,m}(X) =
\bigg( \E \max_{1 \le j \le m} |X_j|_H^p
\bigg)^{1/p}
= 1.
$$
Applying Theorem \ref{thrandomintro} with $\de=\eps$, we get that
if $m \ge C^p n^{p/2} (\log m)^{2/p^*} / \eps^{2}$, then
 $$
\E
 \sup_{y \in B_E} \bigg| \frac{n}{m} \sum_{j=1}^m |\lan X_j, y \ran|^p
  -  \| y \|_{E}^p  \bigg| \le
\eps.
 $$

\medskip
Now, we apply Theorem \ref{thrandomintro} for $K=B_H$
which clearly has modulus of convexity of power
type 2 (i.e. satisfies  inequality (\ref{Clark}) for $q=2$).
In that case, $D = 1$, and
$$\sup_{|y|_H \le 1}
 \E|\lan X, y \ran|^2  = 1 / n
 \quad
 \hbox{ and }
 \quad
\gr_{2,m}(X) =
\bigg( \E \max_{1 \le j \le m} |X_j|_2^2
\bigg)^{1/2}
=1.
$$
Applying Theorem \ref{thrandomintro} for $q=p=2$ and $\de=\eps$,
we get that if $m \ge C^2 n \log m / \eps^2$,
$$
 \E
 \sup_{|y|_{H} \le 1} \bigg| \frac{n}{m} \sum_{j=1}^m |\lan X_j, y
 \ran|^2
  -   |y|_{H}^{2} \bigg| \le \eps.
 $$
Choosing the smallest integer $m$ such that, for a new constant
$\tilde C$,
$$
m \ge \frac{\tilde{C}^p}{\eps^{2}} \ n^{p/2} \
(\log n / \eps^{4/p})^{2/p^*}
$$
we get by Chebychev's inequality that there exist $m$ vectors
 $x_{1}, \ldots, x_{m}$
 of Euclidean norm 1 such that for all $y \in E,$
 $$
 \bigg| \frac{n}{m} \sum_{j=1}^m |\lan x_j, y \ran|^p
  -    \| y \|_{E}^p  \bigg| \le
\eps \|y\|_E^p
 $$
 and
 $$
 \bigg| \frac{n}{m} \sum_{j=1}^m |\lan x_j, y
 \ran|^2
  -   |y|_{H}^{2} \bigg| \le
\eps  |y|_{H}^{2}
 $$
 which gives the desired result.
 \hfill $\Box$
%
%
\section{Isotropic log-concave vectors in $\R^n$}\label{seclogconcave}

We investigate the case of $X$ being an isotropic log-concave
vector in $\reel^n$ (or also a vector uniformly distributed in an
isotropic convex body). Let us recall some definitions and
classical facts about $\log$-concave measures. A probability
measure $\mu$ on $\R^n$ is said to be $\log$-concave if for every
compact sets $A, B$, and every $ \la \in [0,1]$, $\mu(\la A +
(1-\la) B) \ge \mu(A)^{\la} \mu(B)^{1-\la}$. There is always a
Euclidean structure $\lan  \cdot, \cdot \ran $ on $\R^n$ for which
this measure is isotropic, i.e. for every $y\in \R^n$,
$$
\E \lan X,y \ran^2 =
\int_{\R^n} \lan x,y \ran^2 d\mu(x) =  | y |_2^2.
$$
A particular case of a $\log$-concave probability measure is the
normalized
 uniform (Lebesgue)
  measure on a convex body.
 Borell's inequality \cite{Borell} (see also \cite{MilSch, MilPaj}) implies that the linear functionals
$x \mapsto \lan x,y \ran$ satisfy  Khintchine type inequalities
with respect to $\log$-concave probability measures. Namely, if $p
\ge 2$, then for every $y \in \R^n$,
\begin{eqnarray}\label{KKlog}
\left( \E \lan X,y \ran^2 \right)^{1/2} \le \left( \E |\lan X,y
\ran|^p \right)^{1/p} \le C p \left( \E \lan X,y \ran^2
\right)^{1/2},
\end{eqnarray}
or in other words
$$
\| \lan \cdot, y \ran \|_{\psi_1} \le C \left( \E \lan X,y \ran^2 \right)^{1/2}.
$$
We have stated in $(\ref{simplerMs})$ that it is easy to deduce
some information about the parameter $\gr_{p,m}(X)$ from the
behavior of the moment $M_s$ of order $s= \max(p, \log m)$ of the
Euclidean norm of the random vector $X$. These moments were
studied for a random vector uniformly distributed in an isotropic
1-unconditional convex body in \cite{BN}, and for a vector
uniformly distributed in the unit ball of a Schatten trace class
in \cite{GP}, where it was proved that  when $s \le c \sqrt n$,
$M_s$ is of the same order as $M_2$ (up to constant not depending
on $s$). Very recently, Paouris \cite{Paouris05} proved that the
same statement is valid for any $\log$-concave isotropic random
vector in $\R^n$. We state precisely his result.

\begin{thm}  [\cite{Paouris05}]%
{There exist  constants $c,C > 0$ such that for any $\log$-concave
isotropic random vector $X$ in $\R^n$, for any $p \le c \, \sqrt
n$,
$$
\left( \E |X|_2^p \right)^{1/p} \le C \left( \E |X|_2^2 \right)^{1/2}.
$$
}
\end{thm}

\noindent
From this sharp estimate, we will deduce the following
\begin{lemma}\label{max}
Let $X$ be an isotropic log-concave random vector  in $\R^n$ and
let $(X_{j})_{1 \le j \le m}$ be  independent copies of $X$. If $m
\le e^{c \, \sqrt n}$, then for any $p \ge 2$
$$
\gr_{p,m}(X) =
\bigg(\E  \max_{1 \le j \le m} |X_{j}|_{2}^p\bigg)^{1/p}
\le
\, \left\{\begin{array}{l}
  C \sqrt n  \hbox{ if } p \le \log m\\
 C \ p \ \sqrt n  \hbox{ if } p \ge \log m
\end{array}
\right.
$$
\end{lemma}
{\bf Proof.} Since $X$ is isotropic, and for every $y \in \R^n$,
$\E \lan X, y \ran^2 = |y|_2^2$, we get $\E |X|_{2}^2 = n$. By
Borell's inequality \cite{Borell}, $ \forall q \ge 2$, $(\E
|X|_2^q)^{1/q} \le C q \sqrt n$. Therefore if $p \ge \log m$,
$$
\bigg(\E  \max_{1 \le j \le m} |X_{j}|_{2}^p\bigg)^{1/p}
\le
\bigg(\E  \sum_{1 \le j \le m} |X_{j}|_{2}^p\bigg)^{1/p}
\le C p m^{1/p} \sqrt n \le C p \sqrt n.
$$
If $p \le \log m$, by
$(\ref{simplerMs})$
$$
\bigg(\E  \max_{1 \le j \le m} |X_{j}|_{2}^p\bigg)^{1/p}
\le
e \left( \E |X|_2^{\log m} \right)^{1/ \log m}.
$$
Since $m \le e^{c \, \sqrt n}$, $\log m \le c \, \sqrt n$, the
Theorem of Paouris implies
$$
\left( \E |X|_2^{\log m} \right)^{1/ \log m} \le C \sqrt n,
$$
which concludes the proof of the Lemma.
\hfill $\Box$
\begin{corollary}\label{psi1maxlog}
Let $X$ be an isotropic log-concave random vector in $\R^n$, and
let $(X_{j})_{1 \le j \le m}$ be independent copies of $X$. Then
for every $m \le e^{c \, \sqrt n}$
$$
\| \, \max_{1 \le j \le m} |X_j|_2 \ \|_{\psi_1} \le C \sqrt{n}.
$$
\end{corollary}
{\bf Proof.} By Lemma \ref{max}, we know that
$$
\forall r \ge 2, \,
\left( \E \max_{1 \le j \le m} |X_j|_2^r \right)^{1/r} \le C r \sqrt{n}
$$
which proves the claimed estimate for the $\psi_1$-norm.
$\hfill \Box$

\smallskip
\noindent {\bf Remark.} Recall that for a random isotropic
log-concave vector, Borell's inequality implies that
$$
\| \ |X|_2 \ \|_{\psi_1} \le C \sqrt n.
$$
Therefore, a direct application of Lemma \ref{psia} is not enough
to obtain the desired estimate.
%
%

\smallskip
We are now able to give a proof of Theorem \ref{logconcaveintro}.
It is based on the  estimates of $\gr_{p,m}(X)$ proved above.
%
%

\smallskip
\noindent {\bf Proof of Theorem \ref{logconcaveintro}.} Let $\eps
\in (0,1)$ and $p \ge 2$ and set $n_0(\eps, p) = c_p +
\eps^{-4/p}$ where $c_p$ depends only on $p$. For any $n \ge
n_0(\eps, p)$, for any log-concave isotropic random vector $X$ in
$\R^n$, set
$$
V_p =  \sup_{y \in  B_2^n} \bigg| \frac{1}{m} \sum_{j=1}^m |\lan X_j, y \ran|^p
  \  - \  \E |\lan X, y \ran |^p \bigg|
$$
where $X_1, \ldots, X_m$ are independent copies of $X$. Assume
that $p \le \log m$ and $m \le e^{c \, \sqrt n}$ then by Lemma
\ref{max}, we know that
$$
\gr_{p,m}(X)^p \le c_1^p n^{p/2}.
$$
We shall use Theorem \ref{thrandomintro} with $K = B_2^n$ which is
uniformly convex of power type 2 with constant 1 and for which
$D=1$. By  (\ref{KKlog}),
$$
1 \le \sup_{y \in B_2^n} \E | \lan X, y
\ran |^p \le p^p,
$$
therefore Theorem \ref{thrandomintro} implies
that for every $\delta \in (0,1)$, satisfying \newline
 $ C^p n^{p/2} (\log m) \le \delta^{2} m$, we have
$$
\E
 \sup_{y \in B_2^n} \bigg| \frac{1}{m} \sum_{j=1}^m |\lan X_j, y \ran|^p
  \  - \  \E |\lan X, y \ran|^p \bigg| \le 2 \delta p^{p}.
$$
By taking  $\delta$ such that $2 \delta p^{p} = \eps$, we deduce that
$$
\hbox{if }
m \ge C_p' \ \eps^{-2} n^{p/2} \log (n \eps^{-4/p})
\hbox{ then }
\E V_p \le \eps.
$$
Since $n \ge n_0(\eps,p)$, it is easy to see that if
$$
m = \lfloor C_p \ \eps^{-2} n^{p/2} \log n \rfloor,
$$
then $m \ge C_p' \ \eps^{-2} n^{p/2} \log (n \eps^{-4/p})$, $m \le e^{c \sqrt n}$ and $p \le \log m$
which allows us to use the estimate $\E V_p \le \eps$.

To get a deviation inequality for $V_p$, we will apply a  result similar to
Theorem \ref{concentration}. We know by Corollary \ref{psi1maxlog}
that
$$
\| \, \max_{1 \le j \le m} |X|_2 \ \|_{\psi_1} \le C \sqrt{n}.
$$
Following the proof of Theorem \ref{concentration} and replacing inequality
$(\ref{replace})$ by the previous estimate, we easily see that
$$
\| V_p \|_{\psi_{1/p}}
\le
C_p \left(
\E V_p + \frac{2 n^{p/2}}{m}
\right).
$$
Since  $ \lfloor C_p \ \eps^{-2} n^{p/2}  \ \log n \rfloor = m$ then
$$
\E V_p \le \eps \hbox{ and } 2 n ^{p/2}/m \le \eps
$$
 and we deduce  from the
 Chebychev inequality that for any $t > 0$,
 $$
 \Prob ( V_p \ge t ) \le C \exp (- \left( t / C_p' \eps \right)^{1/p}).
 $$
 Therefore, for any $t \ge \eps$,
 with probability greater than $1 - C \exp (- \left( t / C_p' \eps \right)^{1/p})$, $V_p \le t$ which means that
 $$
 \forall y \in \R^n,
 \bigg| \frac{1}{m} \sum_{j=1}^m |\lan X_j, y \ran|^p
  \  - \  \E |\lan X, y \ran|^p \bigg|
\le t |y|_2^p.
 $$
 Since $|y|_2 = (\E \lan X, y \ran^2)^{1/2} \le (\E |\lan X, y \ran |^p)^{1/p}$, we get the claimed result of Theorem
 \ref{logconcaveintro}.
 \hfill $\Box$

 \smallskip
 \noindent
 {\bf Remark.}
Since by Borell's inequality (\ref{KKlog}), for any $ y \in \R^n$,
$$
|y|_2 = \left( \E \lan X,y \ran^2 \right)^{1/2}
\le
\left( \E |\lan X,y \ran|^p \right)^{1/p}
\le
 C p \left( \E \lan X,y \ran^2 \right)^{1/2} = C p |y|_2,
$$
it is clear that Theorem \ref{logconcaveintro} improves the results of Giannopoulos and Milman
\cite{GiaMil}.
%
%
%
\section{When the linear functionals associated
to the random vector $X$ satisfy a $\psi_2$ condition}

Let start this section considering the case when $X$ is a Gaussian
vector in $\R^n$. Let $X_j, \ j= 1, \ldots, m,$ be independent
copies of $X$. For $t \in \R^m$ denote by $X_{t,y}$ the Gaussian
random variable
$$
X_{t,y} = \sum_{j=1}^m t_{j} \lan X_j, y \ran.
$$
Observe that if $p^*$ denotes the conjugate of $p$, then
$$
\sup_{t \in B_{p^*}^m} X_{t,y} =
\bigg(\sum_{j=1}^m |\lan X_j, y \ran|^p \bigg)^{1/p}.
$$
Let $Z$ and $Y$ be  Gaussian vectors in $\R^m$ and $\R^n$
respectively. Using Gordon's inequalities \cite{Gordon}, it is
easy to show that whenever $\E |Z|_p \ge \eps^{-1} \E |Y|_2$ (i.e.
for a universal constant $c$, $m \ge c^p p^{p/2} \eps^{-p}
n^{p/2}$)
$$
\begin{array}{rcl}
\E |Z|_p - \E |Y|_2 \le
\dis
\E \inf_{y \in S^{n-1}}
\bigg(\sum_{j=1}^m |\lan X_j, y \ran|^p \bigg)^{1/p}
& \le &
\\
\le \dis \E \sup_{y \in S^{n-1}} \bigg(\sum_{j=1}^m |\lan X_j, y
\ran|^p \bigg)^{1/p} & \le & \E |Z|_p + \E |Y|_2,
\end{array}
$$
where $(\E |Z|_p + \E |Y|_2) / (\E |Z|_p - \E |Y|_2) \le
(1+\eps)/(1-\eps)$. It is therefore possible to get (with high
probability with respect to the dimension $n$, see \cite{Gordon2})
a family of $m$ random vectors $X_1, \ldots, X_m$ such that for
every $y \in \R^n$,
$$
A \ |y|_2 \le \bigg(\frac{1}{m} \sum_{j=1}^m |\lan X_j, y \ran|^p
\bigg)^{1/p} \le A \ \frac{1 + \eps}{1 - \eps} \ |y|_2.
$$
This argument significantly improves the bound on $m$ in Theorem
\ref{logconcaveintro} for Gaussian random vectors.

 In this part we will be
interested in isomorphic moment estimates (instead of almost
isometric as in Theorem \ref{logconcaveintro}).  We will be able to
extend the estimate for the Gaussian random vector to random
vector $X$ satisfying the $\psi_2$ condition for linear
functionals $y \mapsto \lan X, y \ran$  with the same dependance
on $m$.

Recall that a random variable $Z$ satisfies the   $\psi_2$
condition if and only if for any $\la \in \R$
$$
\label{randpsi2}
  \E \exp (\la Z) \le 2 \exp( c \la^2 \cdot \norm{Z}_2^2).
$$
We prove the following
\begin{theorem}\label{logconcavepsi2}
Let $X$ be an isotropic random vector in $\R^n$ such that all
functionals $y \mapsto \lan X, y \ran$ satisfy the $\psi_2$
condition. Let  $X_1, \ldots, X_m$ be independent copies of $X$.
Then for every $p \ge 2$ and every  $m \ge  n^{p/2}$
$$
\E \sup_{y \in B_2^n} \bigg( \frac{1}{m} \sum_{j=1}^m |\lan X_j, y
\ran|^p \bigg)^{1/p} \le c \, \sqrt p.
$$
\end{theorem}
Note that the results of Part 3 of \cite{GiaMil} follow
immediately from Theorem \ref{logconcavepsi2}, since the random
vector with independent $\pm 1$ coordinates satisfies the $\psi_2$
condition for scalar products.

{\bf Proof.}
 Since $X$ is isotropic,
 \[
   \norm{\pr{X}{y}}_{\psi_2} \le c \norm{\pr{X}{y}}_{2}
   = c |y|.
 \]
Hence, for any $\la \in \R$
$$
\E \exp \la \lan X,y \ran \le 2 e^{c \la^2 |y|_2^{2}}.
$$
Writing
$$
\Delta = X_{t,y} - X_{t',y'} = \sum_{j=1}^m \left( (t_j - t'_j)
\lan X_j,y \ran + t'_j \lan X_j, y - y' \ran \right),
$$
it is easy to find a new constant $c \ge 1$ such that
for every  $t, t' \in B_{p^*}^m,$ $y, y' \in B_2^n$ and every $\la
\in \R^{+}$,
$$
\E \exp( \la \Delta ) \le  2 e^{c \la^2 (|t-t'|_{2}^{2} +
|y-y'|_{2}^{2})}.
$$

This means that $\norm{\Delta}_{\psi_2} \le c (|t-t'|_{2}^{2} +
|y-y'|_{2}^{2})^{1/2}$, and so  $X_{t,y}$ is a sub-Gaussian random
process with respect to the distance
\[
 d\big( (t,y) ; (t',y') \big)
= \big(|t-t'|_{2}^{2} + |y-y'|_{2}^{2} \big)^{1/2}.
\]
Let $G_{t,y} = \langle Z, t \rangle + \langle Y, y \rangle$, where
$Z \in \R^m$ and $Y \in \R^n$ are two independent Gaussian
vectors. Then
$$
\big( \E
|G_{t,y} - G_{t', y'}|^{2} \big)^{1/2}
=
d\big( (t,y) ; (t',y') \big)
$$
The natural metric for the random process $X_{t,y}$ is bounded by
the metric of the process  $G_{t,y}$. The Majorizing Measure
theorem of Talagrand \cite{Tal:genericchaining} implies that
\[
  \E \sup_{(t,y) \in V} X_{t,y}
  \le C \sup_{(t,y) \in V} G_{t,y}
\]
for any compact set $V \subset \R^m \times \R^n$. Therefore,
$$
\begin{array}{l}
\dis
\E \sup_{y \in B_2^n} \bigg(\frac{1}{m} \sum_{j=1}^m |\lan X_j, y \ran|^p \bigg)^{1/p}
=
\frac{1}{m^{1/p}} \E \sup_{t \in B_{p^*}^m} \sup_{y \in B_2^n}
 \sum_{j=1}^m t_{j} \lan X_j, y \ran
\\
 \le
 \dis
 \frac{C}{m^{1/p}}
 \E \sup_{t \in B_{p^*}^m} \sup_{y \in B_2^n}
 G_{t,y} =  \frac{C}{m^{1/p}} \big( \E |Z|_p + \E |Y|_2 \big)
 \\
 \le
 \dis C \big( \sqrt p + \frac{\sqrt n}{m^{1/p}} \big).
\end{array}
$$
This proves that
if $ m \ge  n^{p/2}$, then
$$
\E \sup_{y \in B_2^n} \bigg(\frac{1}{m} \sum_{j=1}^m |\lan X_j, y
\ran|^p \bigg)^{1/p} \le c \sqrt p,
$$
as claimed. \hfill $\Box$

\smallskip
\noindent {\bf Remark.} Let $X$ be an isotropic random vector in
$\R^n$ satisfying the $\psi_2$ estimate for the scalar products.
It is not difficult to see, using Corollary 2.7 in \cite{GiaMil},
that if $m \ge C n$, then with probability greater than $3/4$
$$
c_2 \, |y|_2 \le \bigg( \frac{1}{m} \sum_{j=1}^m |\lan X_j, y
\ran|^p \bigg)^{1/p}
$$
for every $y \in \R^n$. Therefore, using Theorem
\ref{logconcavepsi2}, it is easy to deduce that if $m \ge
n^{p/2}$, then with probability greater than $1/2$
$$
   \forall y \in \R^n \quad
   c_2 \ |y|_2
   \le \bigg( \frac{1}{m} \sum_{j=1}^m |\lan X_j, y \ran|^p
       \bigg)^{1/p}
   \le c_1 \, \sqrt p \, |y|_2
$$
with universal constants $c_1, c_2\ge 1$. This generalizes results
of \cite{GiaMil} and gives an  isomorphic version of the result of
Klartag and Mendelson \cite{KM} valid for every $p\ge 2.$

\bigskip
{\bf Acknowledgement.} A part of this work was done
when the first author was visiting University of Missouri (Columbia). We
wish to thank this institution for its hospitality.

\vspace{1cm}

\noindent Olivier Gu\'edon: \\
 Universit\'e Paris 6, Institut de
 Math\'ematiques de Jussieu, Projet Analyse Fonctionnelle,
 4, place Jussieu, 75005 Paris, France;\\
e-mail: \mbox{guedon@ccr.jussieu.fr}

\vspace{0.5cm}

\noindent Mark Rudelson:  \\
   Department of Mathematics,
   University of Missouri,
   Columbia, MO 65211, USA; \\
   e-mail: \mbox{rudelson@math.missouri.edu}

\end{document}